\documentclass{amsart}

\usepackage{amsmath,amsthm,amsfonts,latexsym,amssymb}
\usepackage{setspace}
\usepackage{ifthen}

\ifthenelse{\pdfoutput > 0}%
{\usepackage[pdftex]{graphicx,color}}%
{\usepackage[dvips]{graphicx,color}}

\usepackage[all]{xy}
\CompileMatrices

\DeclareMathOperator{\Img}{im}

\DeclareMathOperator{\delbar}{\overline{\partial}}
\DeclareMathOperator*{\fibersum}{\#}
\newcommand{\Z}{\mathbb{Z}}
\newcommand{\Q}{\mathbb{Q}}
\newcommand{\A}{\mathbb{A}}
\newcommand{\C}{\mathbb{C}}
\newcommand{\R}{\mathbb{R}}
\newcommand{\D}{\mathbb{D}}

\newtheorem{theorem}{Theorem}[section]
\newtheorem{corollary}[theorem]{Corollary}
\newtheorem{lemma}[theorem]{Lemma}
\newtheorem{lemma*}{Lemma}

\usepackage{hyperref}
\hypersetup{%
  bookmarks=false,%
  pdffitwindow=true,%
  pdftitle={Computations of Floer Homology for certain Lagrangian Tori in closed 4-manifolds},%
  pdfauthor={Adam C. Knapp},%
  pdfnewwindow=true,%
  colorlinks=true,%
  citecolor=black,%
  filecolor=black,%
  linkcolor=black,%
  urlcolor=black,%
}

\title[Computations of Floer Homology]%
{Computations of Floer Homology for certain Lagrangian Tori in closed $4$-manifolds}

\author[A. C. Knapp]{Adam C. Knapp}

\address{Department of Mathematics\\ Michigan State University\\ East
  Lansing, MI 48824}

\email{\href{mailto:knappa@math.msu.edu}{knappa@math.msu.edu}}



\copyrightinfo{2007}{Adam C. Knapp}

\thanks{Partially supported by NSF Grant DMS0305818 and DMS0244663 (FRG)}

\begin{document}
\ifthenelse{\pdfoutput > 0}%
{\DeclareGraphicsExtensions{.pdf,.jpg,.png,.mps}}%
{\DeclareGraphicsExtensions{.mps,.eps,.ps,.jpg,.png}}

\bibliographystyle{mrl}

\begin{abstract}
  We compute the Lagrangian Floer cohomology groups of certain tori in
  closed simply connected symplectic $4$-manifolds arising from
  Fintushel - Stern knot/link surgery. These manifolds are usually not
  symplectically aspherical. As a result of the computation we observe
  examples where $HF(L_0)\cong HF(L_1)$ and $L_0$ and $L_1$ are
  smoothly isotopic but $L_0,L_1$ are not symplectically isotopic and
  are distinguished by $HF(L_0,L_1)$.
\end{abstract}

\maketitle

\section{Introduction}

{\par For a pair $L_0,L_1$ of Lagrangian tori in a symplectic
  $4$-manifold, there are several types of isotopy that we may
  consider:
  \begin{itemize}
  \item Smooth isotopy. 
    
  \item Lagrangian isotopy. $L_0$ and $L_1$ are Lagrangian isotopic in
    $X$ if there is an smooth isotopy of $L_0$ to $L_1$ through
    Lagrangians.

  \item Symplectic/Hamiltonian isotopy. These are Lagrangian isotopies
    of $L_0$ to $L_1$ which extend to isotopies of the ambient
    manifold $X$ and preserve the symplectic structure. If we think of
    a smooth isotopy as given by integrating a time dependant vector
    field $\chi_t$ on $X$, then the condition for symplectic isotopy
    is that $\iota(\chi_t)\omega$ is closed and the condition for
    Hamiltonian isotopy is that $\iota(\chi_t)\omega$ is exact.  Of
    course, in the case of a simply connected ambient manifold, all
    symplectic isotopies are Hamiltonian.

  \end{itemize} }

{\par In \cite{math.GT/0311174}, Vidussi shows that there are an
  infinite number of smoothly nonisotopic Lagrangian tori inside
  $E(2)_K$. In \cite{Fintushel2004}, Fintushel and Stern define an
  integer valued smooth invariant of these tori, distinguishing the
  smooth isotopy class of an infinite family. The examples of
  Lagrangian tori we consider here are of the same type found in
  \cite{math.GT/0311174} and \cite{Fintushel2004}. That is, they occur
  in the ``link portion'' of the link surgery manifolds of
  Fintushel-Stern. }

{\par In \cite{MR1470735} Eliashberg and Polterovich give examples of
  Lagrangian tori in $\R^4$ which are Lagrangian isotopic but not
  Hamiltonian isotopic. The examples considered here differ from these
  as they lie in closed manifolds and may be essential. }

{\par In \cite{MR1743463}, Seidel gives an infinite family of smoothly
  isotopic but symplectically nonisotopic Lagrangian spheres in an
  exact symplectic $4$-manifold distinguished by their Lagrangian
  Floer homology. This result was extended in \cite{MR1765826} to
  certain embeddings of these examples into $K3$ and Enriques
  surfaces. Seidel's computation uses the ``Morse-Bott'' spectral
  sequence for clean intersections of Lagrangians from \cite{Poz} to
  compute the Lagrangian Floer homology of the Lagrangian spheres. We
  say $L_0,L_1$ have clean intersection if their intersection is a
  embedded submanifold and $T(L_0\cap L_1)=TL_0\cap TL_1$. The result
  of Po\'{z}niak in \cite{Poz} was a chain homotopy between the Morse
  complex of the clean intersection and the Lagrangian Floer complex
  restricted to a neighborhood of the clean intersection. We proceed
  somewhat similarly here except that where in Seidel's examples the
  vanishing of higher order differentials comes for more or less for
  free from grading considerations, we show vanishing by an
  essentially topological argument. }

{\par I would like to thank Ronald Fintushel, Richard Hind, and Tom
  Parker for helpful conversations. }

\section{Construction}\label{sec:construct}

{\par Vidussi's symplectic version of the Fintushel-Stern link surgery
  can be described as follows: Consider $M_L$, a $3$ manifold obtained
  from zero surgery on a nontrivial fibered $m$ component link $L$ in
  $S^3$ with a fibration $\pi:M_L\to S^1$.  Choose metrics on $M_L$
  and $S^1$ appropriately so that the fibration map $\pi$ is
  harmonic. (Without loss of generality assume that the metric on
  $S^1$ gives it volume $1$.)  Then $S^1\times M_L$ has a symplectic
  form $\omega = d\theta \wedge d\pi + *_3d\pi$.  Here $\theta$ is a
  coordinate on $S^1$ and $*_3 d\pi$ indicates the pullback of $*d\pi$
  in $M_L$ via the projection $S^1\times M_L \to M_L$. The form
  $\omega$ is closed since $\pi$ is harmonic and nondegenerate since
  $L$ is fibered. If $m_i$ are meridians to the components $K_i$ of
  $L$, then $S^1\times m$ is a symplectic torus of square zero. Let
  $X_{i}$ be symplectic 4-manifolds each with a symplectic torus $F_i$
  of square zero and tubular neighborhood $N(F_i)$. Suppose that
  $\pi_1(X_i\setminus N(F_i))=0$, then the symplectic fiber sum
  $$ X_L=S^1\times M_L \fibersum_{\substack{i = 1,\ldots ,m \\
      F_i=S^1\times m_i}} X_i $$ is simply connected and
  symplectic. Symplectic and Lagrangian submanifolds of each $X_i$ and
  of $S^1\times M_L$ which do not intersect the $F_i = S^1\times
  m_i$ remain so under this process. (We also note that on each link
  component the choice of meridian $m_i$ does not matter since
  isotopies of $m_i$ induce deformation equivalences of symplectic
  structures on $X_L$.)}

{\par Let $\gamma$ be a loop on a fiber of $\pi$. Then with the
  specified symplectic form, $L_\gamma=S^1\times \gamma$ is a
  Lagrangian torus in $S^1\times M_K$. When $\gamma$ and the $m_i$ are
  disjoint, $L_\gamma$ is also naturally a Lagrangian torus in $X_L$.
  It is this class of Lagrangian tori which we will be considering
  here. }

{\par There is one more observation of note. In the $3$-manifold
  $M_L$ there is a natural construction of a vector field $\mu$,
  namely the vector field uniquely determined by
  $\iota(\mu)(*d\pi)\equiv 0$ and $\pi^*(dvol_{S^1})(\mu)\equiv 1$. By
  construction, the time $t$ flow of $\mu$ preserves the fibers of
  $\pi$, moving them in the forward monodromy direction.  Thus the
  time 1 flow of $\mu$ on $M_L$ gives the monodromy map when
  restricted to a fiber of $\pi$. If we extend $\mu$ to a vector field
  on $S^1\times M_L$ which we also call $\mu$, then we note that
  $\iota(\mu)\omega$ is a closed, but not exact 1-form. Thus we get a
  1-parameter family of symplectomorphisms $\phi_t$ on $S^1\times M_L$
  which are not Hamiltonian. }

{\par Consider the action of $\phi_t$ on our Lagrangian torus
  $L_\gamma$. When $t\notin \Z$, $\phi_t(L_\gamma)\cap L_\gamma =
  \emptyset$ as $\gamma$ is moved to a disjoint fiber. However, when
  $t\in \Z$ it is possible that $\phi_t(L_\gamma)$ and $L_\gamma$
  intersect. Further, when the monodromy is of finite order, we can
  find a good choice of meridian $m$ so that the symplectic isotopies
  of $L_\gamma$ to its iterates under the monodromy stay away from
  the $S^1\times m_i$. These symplectic isotopies survive as Lagrangian
  isotopies in $X_L$. }

\section{Calculation of Floer Cohomology}

{\par We use the variant of Lagrangian Floer cohomology over the
  universal Novikov ring $\Lambda$ in \cite{Fukaya}. In this theory,
  the construction of the Floer cohomology groups is defined when
  certain obstruction classes vanish. These classes count
  pseudoholomorphic discs with boundary on the Lagrangians. The
  following lemmas show that we are in the situation where these
  classes vanish and serve to compute the homology. }

\begin{lemma} \label{lemma:nodiscs} Suppose $(S^1\times M_L, \omega)$
  is as above and that the link $L$ is nontrivially fibered in the
  sense that the genus of the fiber is at least $1$.  Then $S^1\times
  M_L$ contains no pseudoholomorphic spheres.  Further, suppose that
  $\gamma_i$, $i=0,1$, are loops on a fiber of $\pi$ which meet
  transversely in exactly one point and let $L_i=L_{\gamma_i}$. Then
  all pseudoholomorphic discs in $S^1\times M_L$ with boundary on
  $L_0$ or on $L_1$ are constant and there are no nonconstant Floer
  discs for $L_0,L_1$.
\end{lemma}

{\par Note that we cannot extend this lemma to say that there are no
  pseudoholomorphic representatives of $\pi_2(S^1\times M_L,L_0\cup
  L_1)$. We see such a counterexample in Section~\ref{sec:ex}. }

\begin{proof}[Proof of Lemma \ref{lemma:nodiscs}]
  \label{proof:nodiscs} As was assumed, $M_L$ is a fibration over
  $S^1$ with fiber $\Sigma_g$ of genus $g\geq 1$ and projection
  $\pi$. As $[\gamma_0]\cdot[\gamma_1]=\pm 1$ in the homology of the
  fiber no nonzero multiples of the two may be homologous. Thus they
  represent distinct infinite order elements of $\pi_1(\Sigma_g,
  \gamma_0\cap\gamma_1)$ for which no powers $i,j\neq 0$ give
  $[\gamma_0]^i=[\gamma_1]^j$. By considering the universal abelian
  cover $\Sigma_g\times \R$, we see that $\pi_1(\Sigma_g)$ injects
  into $\pi_1(M_L)$ by the inclusion of a fiber. Then the subgroups
  generated by $\gamma_0,\gamma_1$ intersect trivially in
  $\pi_1(M_L,\gamma_0\cap \gamma_1)$.

  Since the genus of the fiber $g\geq 1$, $\pi_2(M_L)=0$ and thus
  $\pi_2(S^1\times M_L)=0$.  Now consider the exact sequence
  \begin{displaymath}
    \xymatrix{
      0 = \pi_2(S^1\times M_L) \ar[r] & 
      \pi_2(S^1\times M_L,L_i) \ar[r] & 
      \pi_1(S^1\times \gamma_i) \ar[r]^{i} & 
      \pi_1(S^1\times M_L)
    }
  \end{displaymath}
  It follows that $\pi_2(S^1\times M_L,L_i)=\ker(i)=0$ by our
  assumptions on $\gamma_i$. Therefore, the homomorphism $\int \omega$
  on $\pi_2(S^1\times M_L,L_i)$ given by choosing a representative and
  integrating the pullback of $\omega$ over the disc is the zero
  homomorphism.

  Let $\Omega(L_0,L_1)$ denote the space of paths
  $\delta:([0,1],0,1)\to (S^1\times M_L,L_0,L_1)$ and
  $\Omega_0(L_0,L_1)$ be that subset whose members are homotopic to a
  point.  As $L_0\cap L_1$ is connected, $\Omega_0(L_0,L_1)$ is also
  connected. Let $i_0,i_1$ be the inclusions of $L_0$ and $L_1$ into
  $S^1\times M_L$ respectively.  There is an evaluation map
  $p:\Omega_0(L_0,L_1)\to L_0\times L_1$,
  $p(\delta)=(\delta(0),\delta(1))$. This is a Serre fibration whose
  fiber is homotopy equivalent to $\Omega_0(S^1\times M_L,x)$. Thus
  $\pi_k(p^{-1}(\delta_0,\delta_1)) \cong \pi_{k+1}(S^1\times M_L)$
  and there is the exact sequence
  \begin{displaymath}
    \xymatrix{
      {\pi_2(S^1\times M_L)} \ar[r] \ar[d]^{\cong}&
      {\pi_1(\Omega_0(L_0,L_1))} \ar[r]^{p_*} \ar[d]^{=} &
      {\pi_1(L_0)\times \pi_1(L_1)} \ar[r]^{i_{0*}\cdot i_{1*}^{-1}} \ar[d]^{\cong}&
      {\pi_1(S^1\times M_L)} \ar[d]^{\cong}\\
      0 \ar[r] &
      {\pi_1(\Omega_0(L_0,L_1))} \ar[r] &
      {\Z^2\times \Z^2} \ar[r] &
      {\Z \times \pi_1(M_L) }
    }
  \end{displaymath}
  For the last map, the sequence is exact in the sense
  that $$\Img(p_*) = \ker(i_{0*}\cdot i_{1*}^{-1})= \{ (a,b)\in
  \pi_1(L_0)\times \pi_1(L_1)\mid i_{0*}(a)\cdot (i_{1*}(b))^{-1}=
  e\}.$$ Then $\pi_1(\Omega_0(L_0,L_1)) \cong \ker \left(i_{0*}\cdot
    i_{1*}^{-1}\right)$.  Since the $\gamma_i$ are nontrivial and
  nontorsion each $i_{0*},i_{1*}$ is individually injective on
  $\pi_1$. So the kernel of $i_{0*}\cdot i_{1*}^{-1}$ depends only on
  the intersection of the images of $i_{0*}$ and $i_{1*}$ in
  $\Z\times\pi_{1}(M_L)$.  Then as we know the subgroups generated by
  $\gamma_0,\gamma_1$ in $\pi_1(M_L)$ intersect trivially, we see that
  $\pi_1(\Omega_0(L_0,L_1))\cong \Z$.

  Consider $\D^2$ as the unit disc in $\C$ and let
  $\partial_+=\partial \D^2 \cap \{z \in \C\mid \Re z \geq 0\}$ and
  $\partial_-=\partial \D^2 \cap \{z \in \C\mid \Re z \leq 0\}$.
  Consider $\A$ as the annulus $\{z\in \C \mid 1\leq |z| \leq 2 \}$
  with boundary components $\partial_{|z|=1}$ and $\partial_{|z|=2}$.
  We may represent elements of $\pi_1(\Omega_0(L_0,L_1))$ as maps of
  the annulus $(\A,\partial_{|z|=2},\partial_{|z|=1})$ into
  $(S^1\times M_L,L_0,L_1)$. As with $\pi_2(S^1\times M_L,L_i)$ there
  is a homomorphism on $\pi_1(\Omega_0(L_0,L_1))$ which we shall call
  $\int \omega$ which is given by integrating the pullback of $\omega$
  over (in this case) the annulus. We now show $\int \omega$ on
  $\pi_1(\Omega_0(L_0,L_1))$ is the zero homomorphism.

  To see this, consider a certain generator for
  $\pi_1(\Omega_0(L_0,L_1))\cong \Z$ represented by a map
  $u:(\A,\partial_{|z|=2},\partial_{|z|=1})\to (S^1\times
  M_L,L_0,L_1)$ for which $u(\A)\subset L_0\cap L_1$ and
  $\partial_{|z|=1},\partial_{|z|=2}$ both map to $\pm$ the generator
  of $\pi_1(L_0\cap L_1)=\Z$. Clearly $\int_{\A} u^{*}\omega=0$. Then
  $\int \omega \equiv 0$ on $\pi_1(\Omega(L_0,L_1))$.

  If we have topological Floer disc, that is, a map
  $$u:(\D^2,\partial_-,\partial_+)\to (S^1\times M_L,L_0, L_1),$$ then
  the images of $\pm i$ lie in $L_0\cap L_1$. As $L_0\cap L_1\cong
  S^1$ is connected, we can connect the images of $\pm i$ by an arc
  $\gamma:([-\frac{\pi}{2}, \frac{\pi}{2}], \{-\frac{\pi}{2}\},
  \{\frac{\pi}{2}\})\to (L_0\cap L_1,u(-i),u(i))$.  There are, of
  course, two ways of doing this but that will be immaterial.  Then
  define $\tilde{u}:(\A,\partial_{|z|=2},\partial_{|z|=1}) \to
  (S^1\times M_L,L_0,L_1)$ by
  \begin{displaymath}
    \tilde{u}(z) =
    \left\{
      \begin{array}{ll}
        u\circ \phi & \text{if } \Re z < 0 \\
        \gamma(\theta) & \text{if } \Re z \geq 0, 
        z = re^{i\theta}, -\frac{\pi}{2}\leq \theta \leq \frac{\pi}{2}
      \end{array}
    \right.
  \end{displaymath}
  Where $\phi$ is a diffeomorphism from the interior of $D=\A\cap
  \{z\in \C\mid \Re z \leq 0\}$ to the interior of $\D^2$ whose
  extension to the boundary takes $D\cap \partial_{|z|=2}$ to
  $\partial_-$, $D\cap \partial_{|z|=1}$ to $\partial_+$, the part of
  $D$ lying on the positive imaginary axis to $i$ and the part of $D$
  lying on the negative imaginary axis to $-i$.

  Then $\int_{\A} \tilde{u}^* \omega = \int_{\D^2} u^* \omega$. Now,
  $\tilde{u}$ defines an element of $\pi_1(\Omega(L_0,L_1))$ so
  $\int_{\D^2} u^* \omega = \int_{\A} \tilde{u}^* \omega = 0$. Thus as
  all nonconstant pseudoholomorphic curves have positive symplectic
  area, there are no nonconstant Floer discs.

\end{proof}

\begin{lemma} \label{lemma:nodiscs2} Suppose that $(S^1\times
  M_L,\omega)$ and the $L_i$ satisfy the assumptions of
  Lemma~\ref{lemma:nodiscs}. Let $m_i$ be meridians to each component
  $K_i$ of $L$ each away from the $\gamma_i$ and $(X_i,\omega_{X_i})$
  be a collection of symplectic $4$-manifolds each containing an
  embedded symplectic torus $F_i$ of square zero.  Then for any almost
  complex structure on the fiber sum manifold $X_L$ which on each side
  of the fiber sum is sufficiently close to one for which the $F_i$ or
  $S^1\times m_i$ is pseudoholomorphic, all (perturbed)
  pseudoholomorphic discs in $X_L$ with boundary on $L_0$ or on $L_1$
  and all Floer discs for $L_0,L_1$ are constant.

\end{lemma}  

{\par Note that if $J_t$ is a loop of almost complex structures,
  starting at a $J$ as described, which is contained within a small
  neighborhood of $J$, then $J_t$-pseudoholomorphic strips (Floer
  discs) can be considered as solutions to the perturbed
  pseudoholomorphic curve equations. }

\begin{proof}[Proof of Lemma \ref{lemma:nodiscs2}]
  \label{proof:nodiscs2}
  
  In \cite{MR2113018}, Ionel and Parker construct a $6$ dimensional
  symplectic manifold $Z$ with a map to $\D^2$ so that over
  $\lambda\in \D^2\setminus\{0\}$ the fiber is the symplectic sum
  $X_L$ and over $0$ the fiber is the singular manifold $S^1\times M_L
  \bigcup_{\substack{i=1,\ldots, m \\ F_i=S^1\times m_i}} X_i $. Each
  fiber $X_{L,\lambda}$ is canonically symplectomorphic to $S^1\times
  M_L \setminus N(S^1\times m_i) \cup_{i=1}^{m} X_i\setminus N(F_i)$
  away from the fiber sum region.
  
  Further, the almost-complex structure $J_Z$ on $Z$ is chosen so that
  in the singular fiber $X_{L,0}$, $S^1\times m_i$ is a
  pseudoholomorphic torus and so that the restriction of $J_Z$ to each
  fiber is an almost complex structure. (The singular fiber is
  pseudoholomorphic in the sense that each inclusion of $X$ and
  $S^1\times M_L$ is pseudoholomorphic.) 

  Suppose that $\gamma$ is a smooth embedded path in $\D^2$ which
  passes through zero and that $\{\lambda_n\}_{n=0}^{\infty} \to 0$
  with $\lambda_{n} \in \gamma \setminus \{0\}$. Since the $L_i$ are
  disjoint from the fiber sum region, there are Lagrangian
  submanifolds $\tilde{L}_{i}$ in $Z$ for which in every fiber
  $X_{L,\lambda_{n}}$ above $\lambda_n$, $\tilde{L}_{i}\cap
  X_{L,\lambda_n}= L_{i,n}$ which is mapped to $L_i$ under the
  canonical identification.

  Suppose that we have a family of pseudoholomorphic discs $\Sigma_n$
  in $X_{L,\lambda_n}$ that have boundary in $\tilde{L}_{i}$ for each
  $n$ and all have the same homology class in
  $H_2(X_L,\tilde{L}_i;\Z)$ when it is identified with
  $H_2(X_{L,\lambda_n},L_{i,n};\Z)$ via the bundle over $\gamma$. (The
  projection is a bundle map away from $\lambda=0$.)  Then the
  elements of this family each represent a homology class $[\Sigma_n]$
  in $H_2(Z,\tilde{L}_i;\Z)$. As $Z$ is itself symplectic with a
  symplectic form $\omega_Z$ which restricts to each regular fiber as
  the symplectic form $\omega_{X_{\lambda_n}}$ obtained by fiber sum
  parameterized by $\lambda_n$, the integral of
  $\omega_{X_{\lambda_n}}$ over $\Sigma_n$ within $X_{\lambda_n}$ is
  equal to the integral of $\omega_Z$ over $\Sigma_n$. As $\Sigma_n$
  all represent the same homology class, the integrals must all be
  equal. Then since energy for pseudoholomorphic curves is the
  integral of the symplectic form, we have a bound on energy.
  Then by Gromov compactness for pseudoholomorphic curves with
  Lagrangian boundary conditions, these discs converge to a
  pseudoholomorphic curve $\hat{u}$ with image in the singular fiber
  $X_{L,0} \cong S^1\times M_L \bigcup_{\substack{i=1,\ldots, m \\
      F_i=S^1\times m_i}} X_i$ and boundary on $L_{i,0} \cong
  L_i$. The domain $C$ of $\hat{u}$ is then a collection of spheres
  and discs.

  We next state Lemma 3.4 of \cite{MR1954264} within our context.
  Note that here, $\nu$ is a perturbation of the $\delbar_{J}$
  operator.  (We suppress the perturbation elsewhere.)
  \begin{lemma*}[Ionel,Parker]
    Suppose that $C$ is a smooth connected curve and $f:C \to
    S^1\times M_L$ is a $(J,\nu)$-holomorphic map that intersects
    $S^1\times m$ at a point $p=f(z_0)\in S^1\times m$. Then either
    \begin{enumerate}
    \item $f(C) \subset S^1\times m_i$ for some $i$ or
    \item \label{lem:part2} there is an integer $d>0$ and a nonzero
      $a_0\in \C$ so that in local holomorphic coordinates centered at
      $p$
      \begin{displaymath}
        f(z,\bar{z}) = (p^i + O(|z|),a_0 z^d + O(|z|^{d+1}))
      \end{displaymath}
      where $O(|z|^k)$ denotes a function which vanishes to order $k$
      at $z=0$.
    \end{enumerate}
  \end{lemma*}
  We note that no irreducible component of $\hat{u}$ is mapped
  entirely into any $S^1\times m_i$ as there are no nonconstant
  holomorphic maps of spheres into tori and all discs have boundary on
  $L_i$ away from each $S^1\times m_i$.  Part (\ref{lem:part2}) of the
  lemma shows us that on each component of the domain, $\hat{u}$ meets
  the $F_i$ as algebraic curves do. Thus we know that $\hat{u}$
  intersects each $S^1\times m_i$ at a finite number of points and
  that the image of each spherical or disc component of the domain of
  $\hat{u}$ will lie entirely in one of the $X_i$ or in $S^1\times
  M_L$.

  Then our map $\hat{u}$ splits into two pseudoholomorphic maps
  $\hat{u}_1:(C _1,\partial C_1)\to (S^1\times M_L,L_i)$ and
  $\hat{u}_2:C_2\to \coprod_{i} X_i$. Here $C=C_1 \cup C_2$ with $C_1$ a
  collection of spheres and discs and $C_2$ a (possibly empty)
  collection of spheres so that $C_1\cap C_2 \subset
  \hat{u}^{-1}\left(\bigcup_{i=1}^{m} S^1\times m_i \right)$ are nodes
  of $C$.
  
  Thus we have obtained $\hat{u}_1$, a pseudoholomorphic curve in
  $S^1\times M_L$ with boundary on $L_i$.  By our
  Lemma~\ref{lemma:nodiscs} this map must be constant.  Thus the image
  of $\hat{u}_1$ is in $L_i$ disjoint from the fiber sum region. Then
  since the images of $\hat{u}_2$ and $\hat{u}_1$ are connected,
  $\hat{u}_2$ must have empty domain. Therefore, $\hat{u}$ is
  constant.
  
  This implies that for $\lambda$ sufficiently small, all of the discs
  $\Sigma_n$ must have been contained in $S^1\times M_L$ disjoint from
  the fiber sum region. Therefore by Lemma~\ref{lemma:nodiscs}, they
  were constant in the nonsingular fiber sum.

  The proof follows identically if we consider pseudoholomorphic Floer
  discs for $L_0,L_1$.

\end{proof}

\begin{theorem} \label{thm:SxM_L} Let $(S^1\times M_L,\omega)$ and
  $L_0,L_1$ satisfy the conditions of Lemma~\ref{lemma:nodiscs}.  Then
  \begin{displaymath}
    HF_{S^1\times M_L}(L_0,L_0)
    \cong
    HF_{S^1\times M_L}(L_1,L_1)
    \cong
    H^*(T^2)\otimes \Lambda  
  \end{displaymath}
  and 
  \begin{displaymath}
    HF_{S^1\times M_L}(L_0,L_1)
    \cong
    H^*(S^1)\otimes \Lambda.
  \end{displaymath}
\end{theorem}  

\begin{corollary} \label{cor:SxM_L} Under the conditions of
  Theorem~\ref{thm:SxM_L}, the Lagrangian tori $L_i=S^1\times
  \gamma_i$ are not Hamiltonian isotopic in $(S^1\times M_L,\omega)$.
\end{corollary}

\begin{proof}[Proof of Theorem \ref{thm:SxM_L}] \label{proof:SxM_L}

  We begin by computing $HF(L_0,L_0)$. Since we have found in
  Lemma~\ref{lemma:nodiscs} that $\pi_2(S^1\times M_L,L_0)=0$, we have
  that
  \begin{displaymath}
    HF(L_0,L_0) \cong H^*(T^2) \otimes \Lambda  
  \end{displaymath}
  as in Floer's original work \cite{MR965228}. Similarly, $HF(L_1,L_1)
  = H^*(T^2) \otimes \Lambda$.
  
  Now we consider $HF(L_0,L_1)$. As $L_0,L_1$ intersect cleanly,
  Proposition 3.4.6 of \cite{Poz} implies that in some neighborhood
  $N(L_0\cap L_1)$ of $L_0\cap L_1$ the Floer complex and Morse
  complex for some Morse function $f:L_0\cap L_1\to \R$ coincide.
  This allows us to consider a slight modification of the action
  spectral sequence of \cite{Fukaya}.

  The universal Novikov ring $\Lambda$ can be written as the ring of
  formal sums $\sum_i a_i T^{\lambda_i} e^{n_i}$ with
  \begin{enumerate}
  \item $\lambda_i \in \R$ and $n_i\in \Z$
  \item for each $\lambda^* \in \R$, $\#\{i | \lambda_i \leq
    \lambda^*\} < \infty$
  \end{enumerate}
  Here the $T^{\lambda}$ parameter will be used to keep track of the
  action of a pseudoholomorphic disc. The formula $\deg
  T^{\lambda}e^{n} =2n$ determines a grading on $\Lambda$ and we
  denote by $\Lambda^{k}$ the homogeneous degree $k$ part.

  An $\R^+$ filtration on $\Lambda$ is given by
  \begin{displaymath}
    F^{\lambda}\Lambda =
    \left\{
      \sum_i a^i T^{\lambda_i}e^{n_i} \mid \lambda_i \geq \lambda
    \right\}    
  \end{displaymath}
  We can then get a $\Z$ filtration by picking some $\lambda^*\in\R^+$
  and setting $\mathcal{F}^q \Lambda = F^{q\lambda^*} \Lambda$. The
  homogeneous elements of level $q$ are then
  $\text{gr}_q(\mathcal{F}\Lambda) = \mathcal{F}^{q}\Lambda /
  \mathcal{F}^{q+1}\Lambda$.
  
  Then as in Theorem 6.13 of \cite{Fukaya} we have a spectral sequence
  $E^{p,q}_r$ with 
  \begin{displaymath}
    E^{p,q}_2  
    = 
    \bigoplus_k H^{k}(L_0 \cap L_1;\Q) \otimes \text{gr}_q(\Lambda^{(p-k)})
  \end{displaymath}
  Thus $E_2 \cong H^*(S^1) \otimes \Lambda$. We now want to see that
  all higher order differentials vanish. Recall that in
  Lemma~\ref{lemma:nodiscs} we showed that $\omega|_{\pi_2(S^1\times
    M_L,L_0\cup L_1)}\equiv 0$.  Since the Hamiltonian perturbation
  may be chosen to be very small so that the local curves are of area
  less than $\lambda^*$, Lemma~\ref{lemma:nodiscs} shows that we have
  found all of the discs to be counted.  Therefore, there are no
  higher order differentials and
  \begin{displaymath}
    HF(L_0,L_1)= H^*(S^1)\otimes \Lambda
  \end{displaymath}

\end{proof}

\begin{theorem} \label{thm:XL} Suppose $(X_L,\omega)$ and $L_i$
  satisfy the conditions of Lemma~\ref{lemma:nodiscs2} and that for
  each $L_i$, the Maslov class $\mu:\pi_2(X_L,L_i)\to \Z$ takes only
  even values.  Then in the fiber sum manifold $X_L$,
  \begin{displaymath}
    HF_{X_L}(L_0,L_0)
    \cong
    HF_{X_L}(L_1,L_1)
    \cong 
    H^*(T^2)\otimes \Lambda 
  \end{displaymath}
  and
  \begin{displaymath}
   HF_{X_L}(L_0,L_1)
   \cong 
   H^*(S^1)\otimes \Lambda. 
  \end{displaymath}
\end{theorem}  

\begin{corollary} \label{cor:XL} Under the conditions of
  Theorem~\ref{thm:XL}, the Lagrangian tori $L_i=S^1\times \gamma_i$
  are not symplectically isotopic in $X_L$.
\end{corollary}

{\par In section~\ref{sec:ex}, We give examples which are Lagrangian
  isotopic. }

\begin{proof}[Proof of Theorem \ref{thm:XL}]

  This is proved using the same methods as in the previous theorem.
  Lemma~\ref{lemma:nodiscs2} ensures that the obstructions to defining
  a small perturbation Lagrangian Floer cohomology on the $L_i \subset
  X_L$ vanish and we make make the following local (in $J$) computations:
  \begin{itemize}
  \item Using the action spectral sequence for $HF(L_0,L_0)$ we see
    that $E_2 = H^*(T^2)\otimes \Lambda$. By
    Lemma~\ref{lemma:nodiscs2}, we see that all the higher order
    differentials vanish for a small perturbation. In this case, the
    calculation of $HF(L_0,L_0)$ goes through as in the case of
    Theorem \ref{thm:SxM_L} and we find that $HF(L_0,L_0) \cong
    H^*(T^2) \otimes \Lambda$.  Similarly, we get $HF(L_1,L_1) \cong
    H^*(T^2) \otimes \Lambda$.

  \item Using the action spectral sequence for $HF(L_0,L_1)$ we see
    that $E_2 = H^*(S^1)\otimes \Lambda$. By
    Lemma~\ref{lemma:nodiscs2}, we see that all the higher order
    differentials vanish. Therefore, the calculation of $HF(L_0,L_1)$
    goes through as in the case of Theorem \ref{thm:SxM_L} and we find
    that $HF(L_0,L_1) \cong H^*(S^1) \otimes \Lambda$.

  \end{itemize}
  

  Now we show that the computed groups are invariant. Invariance is
  guaranteed as long as in $1$-parameter families, there are no discs
  with boundary on $L_0,L_1$ of index $-1$ which bubble
  off.
  

  Our restriction on the Maslov class ensures that there are no index
  $-1$ discs which appear in the boundary of $1$-parameter
  families. Thus there is a continuation isomorphism and we have a
  well defined invariant of Hamiltonian isotopy in $X_L$.

\end{proof}

\section{Example}\label{sec:ex}

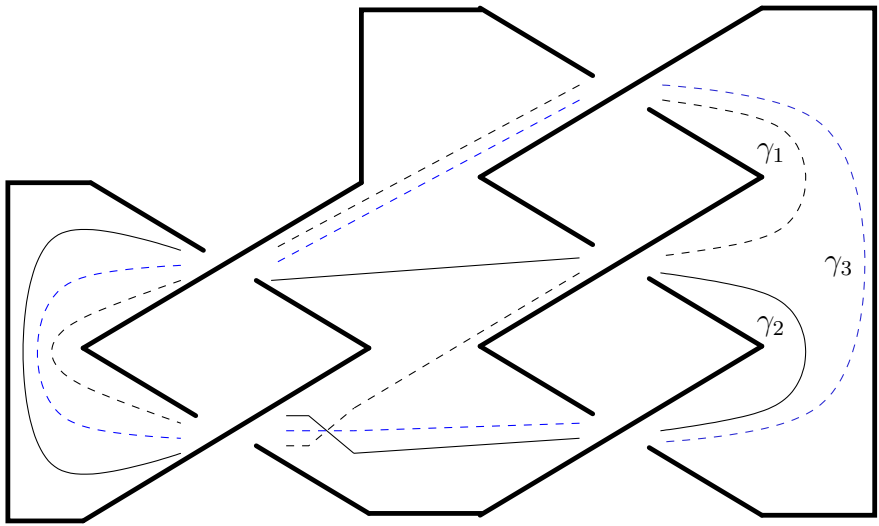
\begin{figure}[htb]
  \centering
  \begin{picture}(0,0)%
    \includegraphics{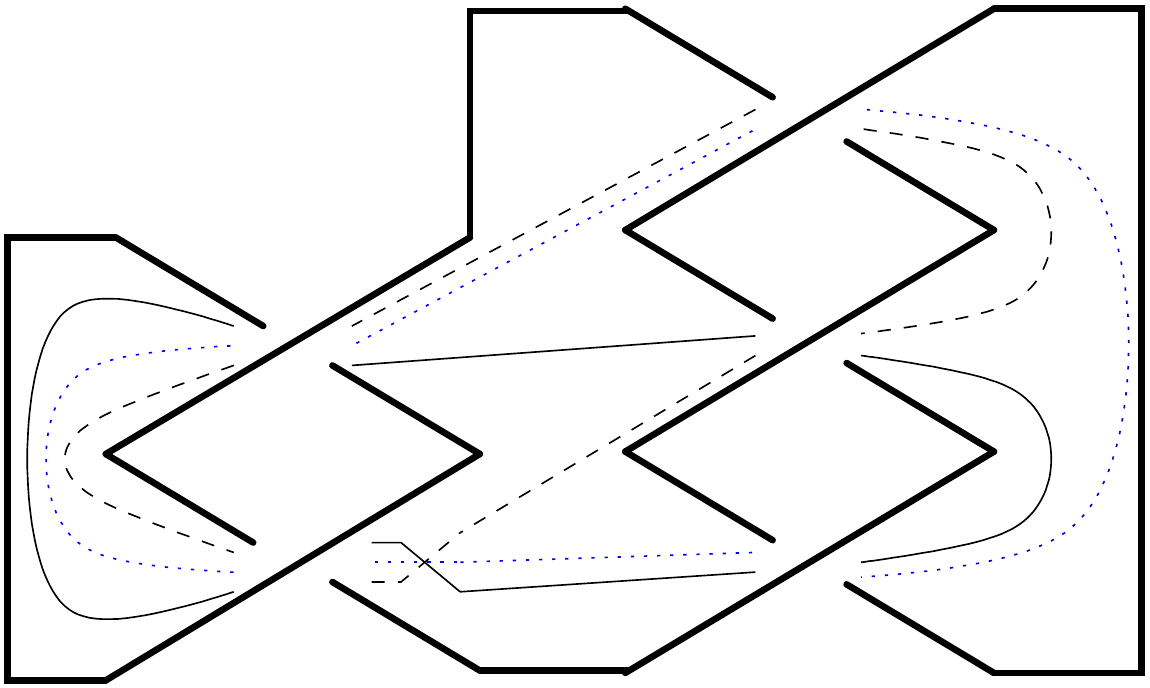}%
  \end{picture}%
  \setlength{\unitlength}{4144sp}%
  \begingroup\makeatletter\ifx\SetFigFontNFSS\undefined%
  \gdef\SetFigFontNFSS#1#2#3#4#5{%
    \reset@font\fontsize{#1}{#2pt}%
    \fontfamily{#3}\fontseries{#4}\fontshape{#5}%
    \selectfont}%
  \fi\endgroup%
  \begin{picture}(5252,3137)(913,-2254)
    \put(5581,-61){\makebox(0,0)[rb]{\smash{{\SetFigFontNFSS{12}{14.4}{\familydefault}{\mddefault}{\updefault}{\color[rgb]{0,0,0}$\gamma_1$}%
          }}}}
    \put(5986,-736){\makebox(0,0)[rb]{\smash{{\SetFigFontNFSS{12}{14.4}{\familydefault}{\mddefault}{\updefault}{\color[rgb]{0,0,0}$\gamma_3$}%
          }}}}
    \put(5581,-1096){\makebox(0,0)[rb]{\smash{{\SetFigFontNFSS{12}{14.4}{\familydefault}{\mddefault}{\updefault}{\color[rgb]{0,0,0}$\gamma_2$}%
          }}}}
  \end{picture}%
  \caption{$\gamma_i$, $i=1,2,3$ intersecting pairwise transversely
    (within the fiber) in single points}
  \label{fig:knot}
\end{figure}

\begin{figure}[p]
  \centering
\begin{picture}(0,0)%
\includegraphics{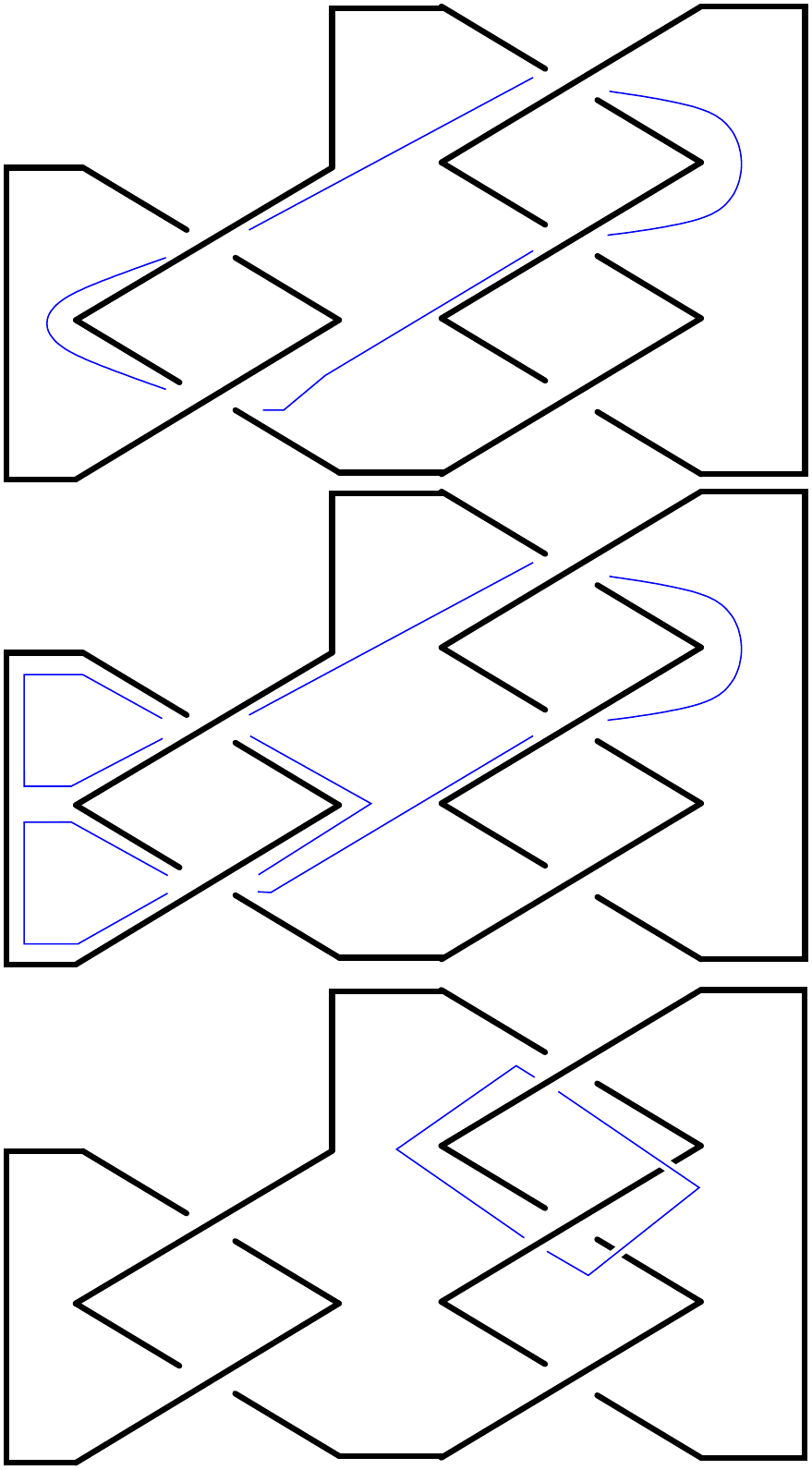}%
\end{picture}%
\setlength{\unitlength}{3947sp}%
\begingroup\makeatletter\ifx\SetFigFontNFSS\undefined%
\gdef\SetFigFontNFSS#1#2#3#4#5{%
  \reset@font\fontsize{#1}{#2pt}%
  \fontfamily{#3}\fontseries{#4}\fontshape{#5}%
  \selectfont}%
\fi\endgroup%
\begin{picture}(4216,7632)(-9,-6787)
\put(729,439){\makebox(0,0)[b]{\smash{{\SetFigFontNFSS{10}{12.0}{\rmdefault}{\mddefault}{\updefault}{\color[rgb]{0,0,0}$(1)$}%
}}}}
\put(841,-2101){\makebox(0,0)[b]{\smash{{\SetFigFontNFSS{10}{12.0}{\rmdefault}{\mddefault}{\updefault}{\color[rgb]{0,0,0}$(2)$}%
}}}}
\put(852,-4645){\makebox(0,0)[b]{\smash{{\SetFigFontNFSS{10}{12.0}{\rmdefault}{\mddefault}{\updefault}{\color[rgb]{0,0,0}$(3)$}%
}}}}
\put(3661,419){\makebox(0,0)[b]{\smash{{\SetFigFontNFSS{10}{12.0}{\familydefault}{\mddefault}{\updefault}{\color[rgb]{0,0,0}$\gamma_1$}%
}}}}
\end{picture}%
  \caption{Isotopy of $\gamma_1$}
  \label{fig:iso-g1}
\end{figure}

\begin{figure}[p]
  \centering
\begin{picture}(0,0)%
\includegraphics{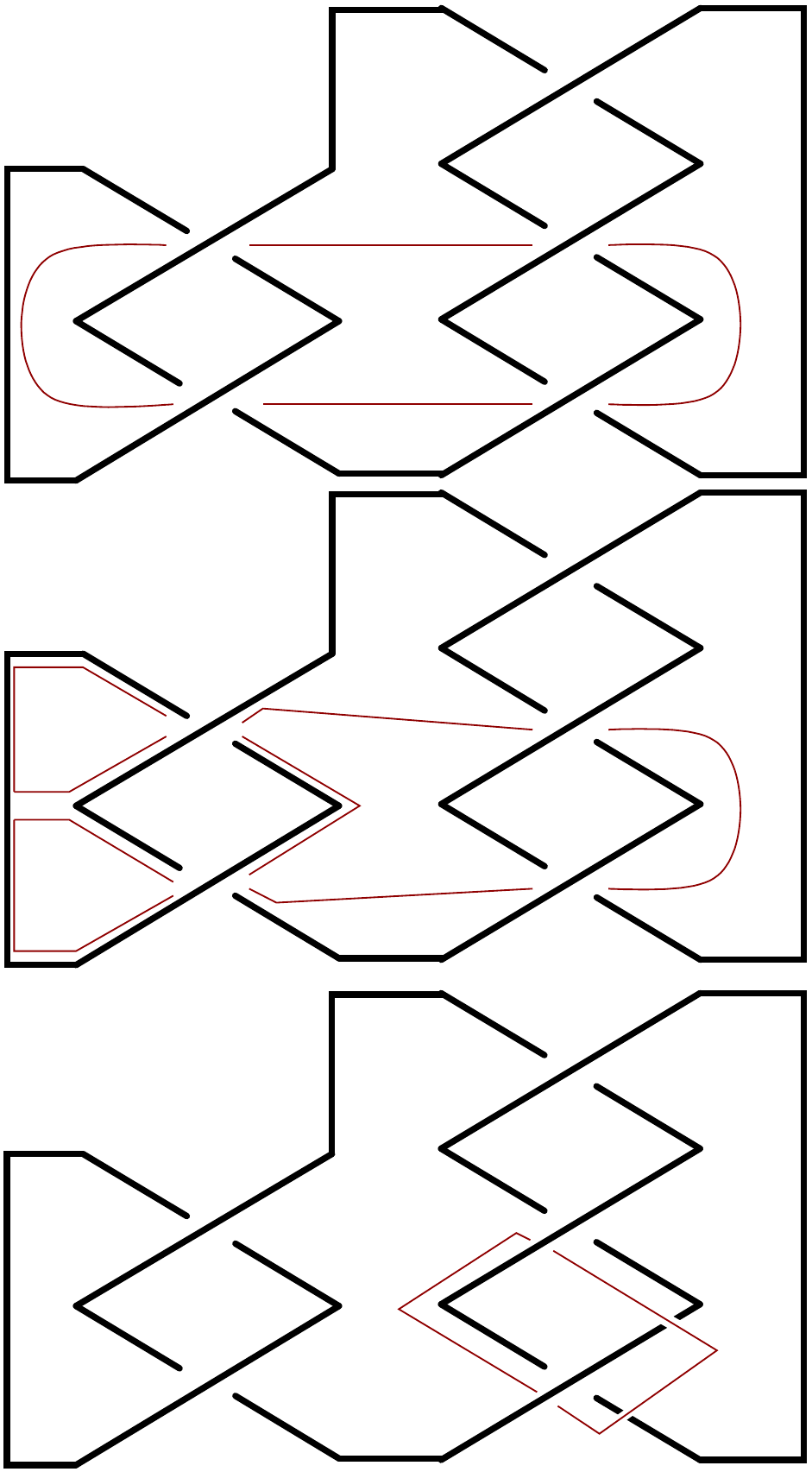}%
\end{picture}%
\setlength{\unitlength}{4641sp}%
\begingroup\makeatletter\ifx\SetFigFontNFSS\undefined%
\gdef\SetFigFontNFSS#1#2#3#4#5{%
  \reset@font\fontsize{#1}{#2pt}%
  \fontfamily{#3}\fontseries{#4}\fontshape{#5}%
  \selectfont}%
\fi\endgroup%
\begin{picture}(3837,6955)(-5,-6107)
\put(991,-1861){\makebox(0,0)[b]{\smash{{\SetFigFontNFSS{10}{12.0}{\rmdefault}{\mddefault}{\updefault}{\color[rgb]{0,0,0}$(2)$}%
}}}}
\put(991,389){\makebox(0,0)[b]{\smash{{\SetFigFontNFSS{10}{12.0}{\rmdefault}{\mddefault}{\updefault}{\color[rgb]{0,0,0}$(1)$}%
}}}}
\put(991,-4336){\makebox(0,0)[b]{\smash{{\SetFigFontNFSS{10}{12.0}{\rmdefault}{\mddefault}{\updefault}{\color[rgb]{0,0,0}$(3)$}%
}}}}
\put(3463,-272){\makebox(0,0)[b]{\smash{{\SetFigFontNFSS{10}{12.0}{\familydefault}{\mddefault}{\updefault}{\color[rgb]{0,0,0}$\gamma_2$}%
}}}}
\end{picture}%
  \caption{Isotopy of $\gamma_2$}
  \label{fig:iso-g2}
\end{figure}
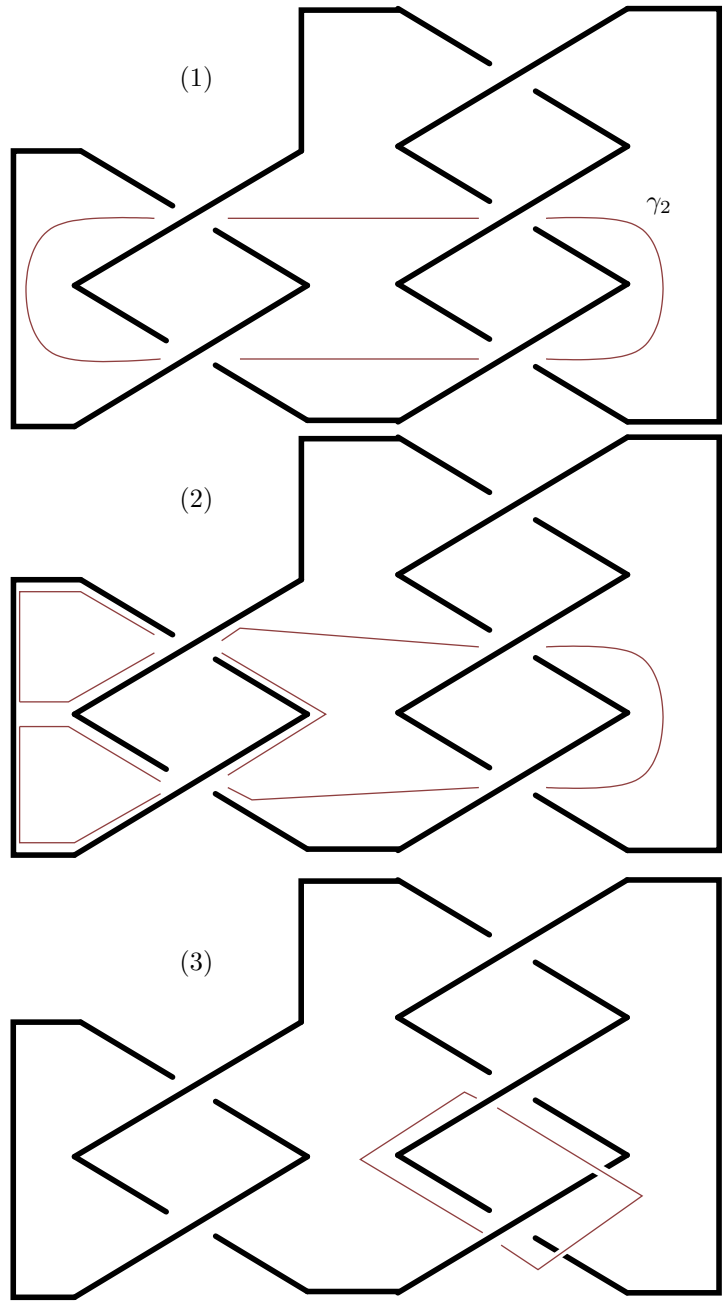

If we let $L$ be the union of the right-hand trefoil knot $K_1$ and
one of its meridians $K_2$, the loops $\gamma_1,\gamma_2,\gamma_3$ in
Figure~\ref{fig:knot} are all {\em freely} smoothly isotopic in $M_L$
and meet transversely pairwise.\footnote{Though they are freely
  isotopic, the proof of Lemma~\ref{lemma:nodiscs} shows that they are
  not equal in $\pi_1$ i.e. fixing a basepoint.} We select the
fibration $\pi:M_L\to S^1$ so that the Seifert surface containing
$\gamma_i$ shown in Figure~\ref{fig:knot} is the fiber. The smooth
isotopies of $\gamma_1$ and $\gamma_2$ to a common curve are shown in
Figures~\ref{fig:iso-g1} and~\ref{fig:iso-g2}. In each of these
figures going from $(1)$ to $(2)$ involves sliding over the
$0$-surgery on the meridian, going from $(2)$ to $(3)$ is an
isotopy. To relate the end results we have the additional move of
``twisting up the corkscrew'' which takes the curves each $(3)$ to the
other. Note that this smooth isotopy is different than the Lagrangian
isotopy which we will mention later.

Since the $\gamma_i$ are smoothly isotopic, the Lagrangian tori
$L_i=S^1\times \gamma_i$ are smoothly isotopic. As Lagrangians do not
have canonical orientations we neglect the orientations of the loops
here.

Addressing the comment made after Lemma~\ref{lemma:nodiscs}, we note
that we can choose the almost complex structure $J$ so that the
Seifert surface $\Sigma$ (a $T^2$) is pseudoholomorphic.  Then
choosing any pair $\gamma_i,\gamma_j$ ($i\neq j$), $\Sigma\setminus
(\gamma_i\cup \gamma_j)$ is a disc. This is however, not a Floer disc
as it does not satisfy the correct boundary conditions.

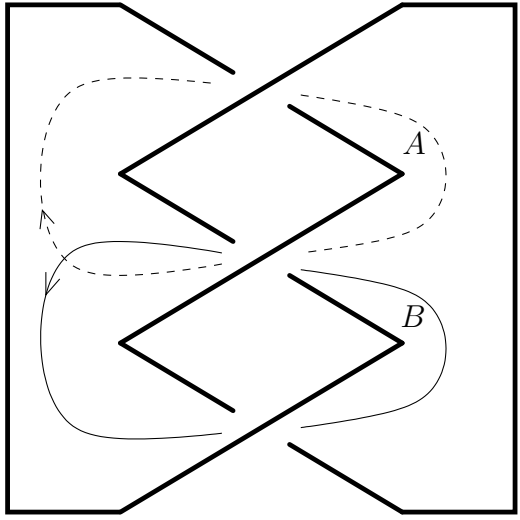
\begin{figure}
  \centering
  \begin{picture}(0,0)%
    \includegraphics{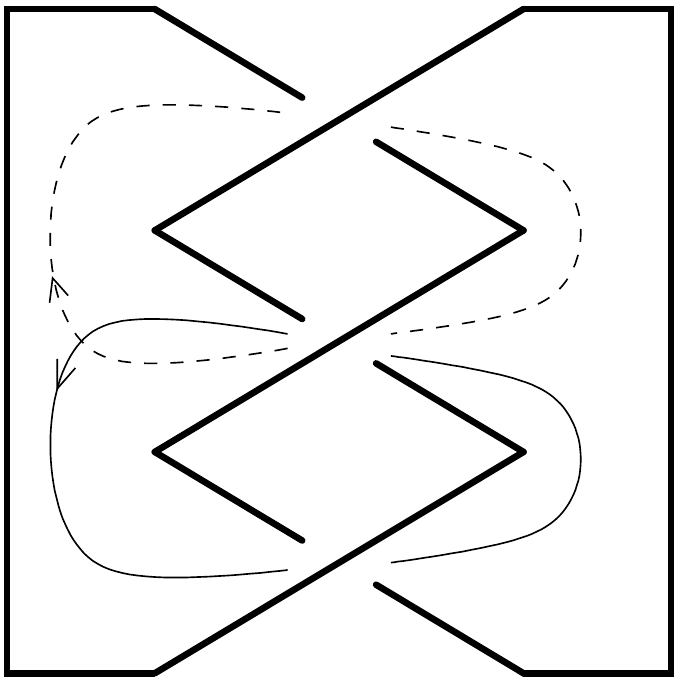}%
  \end{picture}%
  \setlength{\unitlength}{4144sp}%
  \begingroup\makeatletter\ifx\SetFigFontNFSS\undefined%
  \gdef\SetFigFontNFSS#1#2#3#4#5{%
    \reset@font\fontsize{#1}{#2pt}%
    \fontfamily{#3}\fontseries{#4}\fontshape{#5}%
    \selectfont}%
  \fi\endgroup%
  \begin{picture}(3102,3103)(3063,-2220)
    \put(5581,-61){\makebox(0,0)[rb]{\smash{{\SetFigFontNFSS{12}{14.4}{\familydefault}{\mddefault}{\updefault}{\color[rgb]{0,0,0}$A$}%
          }}}}
    \put(5581,-1096){\makebox(0,0)[rb]{\smash{{\SetFigFontNFSS{12}{14.4}{\familydefault}{\mddefault}{\updefault}{\color[rgb]{0,0,0}$B$}%
          }}}}
  \end{picture}%
  \caption{Basis for monodromy on trefoil}
  \label{fig:trefoil}
\end{figure}

For the left and right handed trefoils, the monodromy is of order 6
and in the basis $A,B$ shown in Figure~\ref{fig:trefoil} is given by
the matrix
\begin{displaymath}
  \left[
  \begin{array}{rr}
    1 & 1 \\
    -1 & 0
  \end{array}
  \right]
\end{displaymath}
On the (positive/negative) Hopf link the monodromy is a
(positive/negative) Dehn twist about a curve parallel to the
components. The connect sum of fibered links is fibered with monodromy
which splits around the connect sum region.

From this computation of monodromy, we see that $\gamma_2$ is $\pm$
the 2nd and 5th image of $\gamma_1$ under the monodromy map and that
$\gamma_3$ is $\pm$ the 1st and 4th image of $\gamma_1$. Finally,
$\gamma_1$ is sent to $-\gamma_1$ under the third iteration of the
monodromy.  Thus, as the monodromy gives a symplectic isotopy
(c.f. Section~\ref{sec:construct}) in $S^1\times M_L$ the tori
$L_{1},L_{2},L_{3}$ are all symplectically isotopic.  In $X_L$ they
are all Lagrangian isotopic as the fiber sum is taken away from the
isotopy. Despite the existence of these symplectic isotopies in
$S^1\times M_L$, Theorem~\ref{thm:SxM_L} shows that the
$L_i=S^1\times\gamma_i$ are not Hamiltonian isotopic there.


Now we consider how Theorem~\ref{thm:XL} applies. That all the
hypotheses of the theorem are satisfied, except that on the Maslov
class, is clear. With the following lemma we see that the remaining
condition is satisfied. Then Theorem~\ref{thm:XL} shows that the $L_i$
are not symplectically isotopic in $X_L$.

\begin{lemma}
  For the Lagrangian tori $L_i$, we may choose $X_i=E(1)$, and the
  particular identification of $F_i$ and $S^1\times m_i$ so that
  $\mu_{L_i}:\pi_2(X_L,L_i)\to \Z$ is even.
\end{lemma}

\begin{proof}

  As $L$ is a two component link with odd linking number, $X_L$ is a
  homotopy $E(2)$ and thus is spin. See~\cite{MR1650308}. In fact
  $E(2)_L$ is $E(2)_{\text{Trefoil}}$. We shall write
  $E(2)_L:=X_L$. For such a $4$-manifold, the first Chern class is an
  even multiple of the fiber.

  Note that $\mu_{L_i}$ factors through $\pi_2(E(2)_L,L_i)\to
  H_2(E(2)_L,L_i)$. As $H_1(E(2)_L)=0$, the Meyer-Vietoris sequence
  gives that the group $H_2(E(2)_L,L_i)$ is generated by elements of
  $H_2(E(2)_L)/\left<L_i\right>$ and relative classes with boundary
  spanning $H_1(L_i)$.  The Maslov index of a class $\beta$,
  $\mu_{L_i}(\beta)$, $\beta\in H_2(E(2)_K,L_i)$, will change by an
  even amount whenever an element of $H_2(E(2)_L)/\left<L_i\right>$ is
  added as $c_1(E(2)_L)$ is divisible by $2$. Thus if we can find a
  pair of relative discs whose boundaries generate $H_1(L_i)$ and
  whose Maslov indices are even, we have shown that $\mu$ is even.

  For $i=1,2$, we will choose the identification of $F_i$ and
  $S^1\times m_i$ so that
  \begin{enumerate}
  \item $\text{pt}\times m_i$ is identified with a vanishing cycle on
    $F_i$ and
  \item $S^1 \times \text{pt}$ is identified with the sum of two
    vanishing cycles on $F_i$ whose boundaries meet once,
    transversely, in $F_i$.
  \end{enumerate}
  Because $\pi_1(E(1)\setminus F_i)=1$, we may select an elliptic
  fibration on $E(1)$ with nodal fibers having vanishing cycles $a$
  and $b$ where $a,b$ generate $\pi_1(F_i)$. With the decomposition
  $F_i = a\times b$, identify $a$ with $\text{pt}\times m$ and $a+b$
  with $S^1\times \text{pt}$. This gives us the desired identification
  of $F_i$ and $S^1\times m_i$.

  \begin{figure}[htb]
    \centering
    \includegraphics{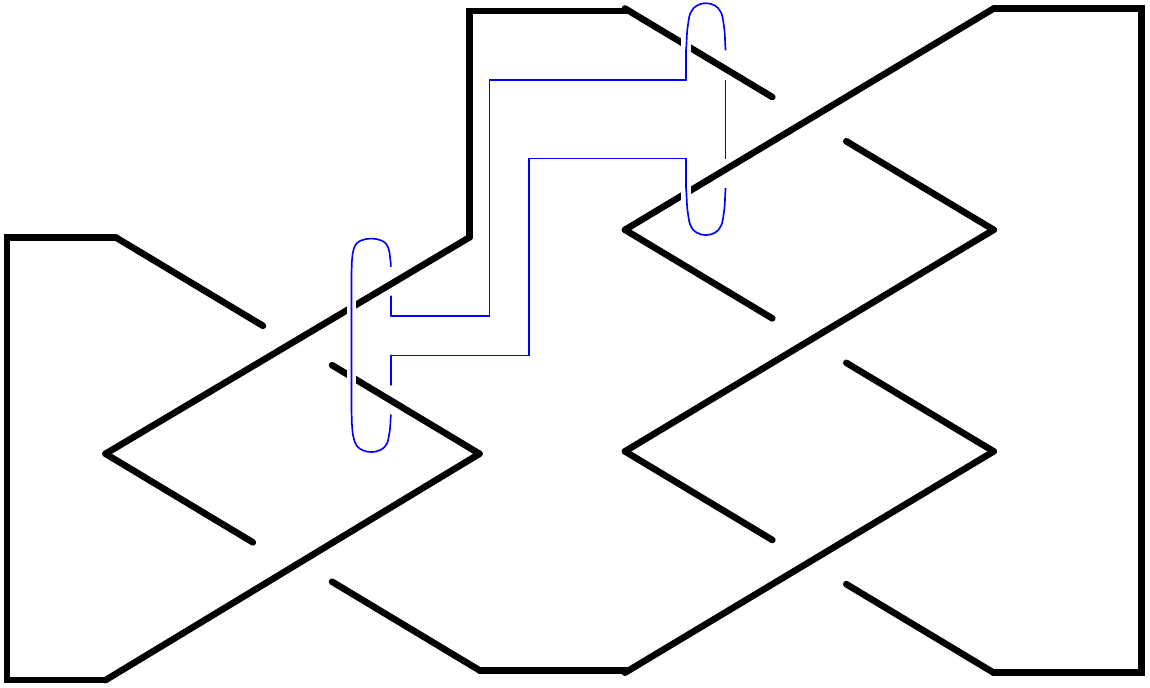}
    \caption{An isotope of $\gamma_1$ bounding meridians to $K_1$ and
      $K_2$}
    \label{fig:g1-cap}
  \end{figure}
  
  Each of the $\gamma_i$ bounds a four times punctured disc
  $D_{0,\gamma_i}$ in $M_L$ where three punctures are meridians of
  $K_1$ and one is a meridian to $K_2$. See Figure~\ref{fig:g1-cap}.
  By our choice of fiber sum gluing, each meridian of $K_1$ bounds a
  vanishing disc on its $E(1)$ side of the fiber sum. Similarly, the
  meridian to $K_2$ bounds a vanishing disc on its $E(1)$ side of the
  fiber sum. Take three copies of the vanishing disc
  $D_{1},D_{2},D_{3}$ from the $E(1)$ fiber summed to $S^1\times m_1$
  and one vanishing disc $D_{4}$ of from the $E(1)$ fiber summed to
  $S^1\times m_2$ and form $D_{\gamma_i}=D_{0,\gamma_i}\cup D_{1}\cup
  D_{2} \cup D_{3} \cup D_{4}$.

  Each $\gamma_i$'s Lagrangian framing relative to that induced by
  trivializing over $D_{\gamma_i}$ is $-2$ and is given by a pushoff
  in the Seifert surface. See~\cite{Fintushel2004}. The framing coming
  from this disc is then $-1-1-1-1-(-2)=-2$ and gives us
  $\mu_{\gamma_i}(D_{\gamma_i})=-2$.

  By our choice of gluing, the $S^1\times\text{pt}\subset L_i$ is
  bounded by a pair of vanishing discs. This pair of vanishing discs
  intersects at one point on $F$ and so can be smoothed to a disc
  $D_{S^1\times\text{pt}}$ with relative framing $-2$. For this loop
  the framing defect from the torus is $0$ (given by pushoff in the
  monodromy direction) and so $\mu_{L_i}(D_{S^1\times\text{pt}})=-2$.

  Thus we have found a basis on which $\mu$ is even, hence $\mu$ is
  even.
  
\end{proof}

This example generalizes to similar links with $K_1=T_{2,2n+1}$ where
we find many Lagrangian but not symplectically isotopic tori.

\end{document}